\documentclass[11pt]{article}
\usepackage{amsmath,amssymb}
\newtheorem{theo}{Theorem}[section]

\newtheorem{lemma}[theo]{Lemma}
\newtheorem{coro}[theo]{Corollary}

\newcommand{\ignore}[1]{}
\def\square{\vrule height6pt width7pt depth1pt}
\def\endpf{\hfill\square\bigskip}

\begin{document}
\title{Packing $4$-cycles in Eulerian and bipartite Eulerian tournaments with an application
to distances in interchange graphs}
\author{Raphael Yuster
\thanks{
e-mail: raphy@research.haifa.ac.il \qquad
World Wide Web: http:$\backslash\backslash$research.haifa.ac.il$\backslash$\symbol{126}raphy}
\\ Department of Mathematics\\ University of
Haifa at Oranim\\ Tivon 36006, Israel}

\date{} 

\maketitle
\setcounter{page}{1}
\begin{abstract}
We prove that every Eulerian orientation of $K_{m,n}$ contains
$\frac{1}{4+\sqrt{8}}mn(1-o(1))$ arc-disjoint directed $4$-cycles, improving earlier lower bounds.
Combined with a probabilistic argument, this result is used to
prove that every regular tournament with $n$ vertices contains
$\frac{1}{8+\sqrt{32}}n^2(1-o(1))$ arc-disjoint directed $4$-cycles.
The result is also used to provide an upper bound for the distance between two antipodal vertices in interchange
graphs.
\end{abstract}

\section{Introduction}
All graphs and digraphs considered here are finite, and contain no parallel edges or anti-parallel arcs.
For standard graph-theoretic terminology the reader is referred to \cite{Bo}.
An {\em Eulerian orientation} of an undirected graph $G$ is an orientation of the edges of $G$
such that each vertex has the same indegree and outdegree in the resulting oriented graph.
It is well known that a graph $G$ has an Eulerian orientation if and only if every vertex is of even
degree. A {\em tournament} is an orientation of a complete graph and a bipartite tournament
is an orientation of a complete bipartite graph. Properties of Eulerian tournaments (also called
{\em regular} tournaments) and
bipartite Eulerian tournaments have been extensively studied in the literature (see, e.g., \cite{Be,Mc,Ro}).

A conjecture of Brualdi and Shen \cite{BrSh} asserts that every bipartite Eulerian tournament
has a decomposition into $C_4$, the directed cycle on four vertices. In other words, any
Eulerian orientation of $K_{m,n}$ contains $mn/4$ arc-disjoint $C_4$'s (clearly $m$ and $n$ must
be even in order to have an Eulerian orientation). This conjecture is still far from being solved,
and, in fact, there is no evidence that it should be true. Even an asymptotic version of this conjecture
is not known. A recent result from \cite{ShYu} states that any Eulerian orientation
of a bipartite graph with vertex class sizes $m$ and $n$ and more than $2mn/3$ arcs has a $C_4$.
This immediately implies that any Eulerian orientation of $K_{m,n}$ contains at least $mn/12$
arc-disjoint $C_4$'s. Although very far from the Brualdi-Shen conjecture, this is currently the best
known published lower bound for this problem. In this paper we considerably improve the lower bound.
\begin{theo}
\label{t1}
Every Eulerian orientation of $K_{m,n}$ contains at least $\frac{1}{4+\sqrt 8}mn(1-o(1))$ 
arc-disjoint copies of $C_4$.
\end{theo}
The $o(1)$ term denotes a function tending to zero as $\min\{m,n\}$ tends to infinity.
Notice that, trivially, any $C_4$-packing of an orientation of $K_{m,n}$ contains at most $mn/4$
elements. Theorem \ref{t2} shows that more than $58$ per cent of the arcs can be packed
with $C_4$'s.

The motivation behind the Brualdi-Shen conjecture stems from the theory of interchange graphs.
These graphs are defined as follows:
Let $R = (r_1,\ldots,r_n)$ and $S = (s_1,\ldots,s_m)$ be non-negative integral vectors
with $\sum r_i=\sum s_j$. Let ${\cal A}(R,S)$ denote the set of all $m \times n$ $\{0,1\}$-matrices with
row sum vector $R$ and column sum vector $S$, and assume that ${\cal A}(R, S) \neq \emptyset$.
This set has been studied extensively (see \cite{Br} for a survey).
An {\em interchange} is a transformation which replaces a $2 \times 2$ identity submatrix of a matrix $A$
with the $2 \times 2$ permutation matrix $P_{i,j}=|i-j|$.
Clearly an interchange (and hence any sequence of interchanges)
does not alter the row and column sum vectors of a matrix and
transforms a matrix in ${\cal A}(R, S)$ into another matrix in ${\cal A}(R, S)$. Ryser \cite{Ry}
proved that for any two matrices in ${\cal A}(R, S)$ there is a sequence of interchanges which
transforms one to the other.
The {\em interchange graph} $G(R,S)$ of ${\cal A}(R,S)$, defined by Brualdi in 1980, is the graph
with all matrices in ${\cal A}(R,S)$ as its vertices, where two matrices are adjacent provided
that one can be obtained from the other by a single interchange.
Brualdi conjectured that the diameter of $G(R,S)$,
denoted $d(R,S)$, cannot exceed $mn/4$.
This conjecture is still far from being resolved. The best known bounds for $d(R,S)$ use
the following reduction of Walkup. Given two matrices $A, B \in {\cal A}(R, S)$ consider the $m \times n$
matrix $A-B$. This matrix contains the elements $\{-1,0,1\}$ and the sum of each row and column is zero.
Thus, there is a one-to-one correspondence between such matrices and Eulerian orientations of
a bipartite graph with vertex class sizes $m$ and $n$.
Walkup \cite{Wa} proved that the distance between $A$ and $B$ in $G(R,S)$, denoted $i(A,B)$ satisfies
\begin{equation}
\label{e1}
i(A,B)=\frac{d(A,B)}{2}-q(A,B)
\end{equation}
where $d(A,B)$ is the number of nonzero entries in $A-B$ and $q(A,B)$ is the maximum number
of arc-disjoint cycles in a cycle decomposition of the bipartite Eulerian digraph corresponding to $A-B$.
The currently best known upper bound for $d(R,S)$ obtained using Walkup's reduction
is given in \cite{ShYu} yielding $d(R,S) \leq \frac{5}{12}mn$.
We call two matrices $A, B \in {\cal A}(R, S)$ {\em antipodal} if $B$ is the complement of $A$,
or, in other words, if $A+B$ is the all-one matrix. Notice that a necessary conditions for $A$ and $B$
to be antipodal is that all coordinates of $R$ are $m/2$ and all coordinates of $S$ are $n/2$.
Trivially, the distance between two antipodal vertices in $G(R,S)$ is at least $mn/4$, and hence such
vertices show that Brualdi's conjecture, if true, is best possible. In particular antipodal vertices are conjectured to
be farthest apart in $G(R,S)$. If $A$ and $B$ are antipodal matrices then $A-B$ corresponds to
an Eulerian orientation of $K_{m,n}$. Observe that using (\ref{e1}) and Theorem \ref{t1} we obtain the following.
\begin{coro}
\label{c1}
Let $A$ and $B$ be two antipodal vertices of $G(R,S)$. Then,
$$
0.25mn \leq i(A,B) \leq \frac{1+\sqrt 2}{4 + \sqrt 8}mn(1+o(1)) \approx 0.354mn(1+o(1)).
$$
\end{coro}

Theorem \ref{t1}, combined with a probabilistic argument, enables us to prove an analogous theorem
for regular tournaments.
\begin{theo}
\label{t2}
Every regular tournament with $n$ vertices contains at least $\frac{1}{8+\sqrt{32}}n^2(1-o(1))$
arc-disjoint copies of $C_4$.
\end{theo}
Notice that, trivially, any $C_4$-packing of an $n$-vertex tournament contains less than $n^2/8$
elements. Theorem \ref{t2} shows that more than $58$ per cent of the arcs can be packed
with $C_4$'s.

In the following section we prove Theorem \ref{t1}. The final section contains the proof of Theorem \ref{t2}.

\section{Proof of Theorem \ref{t1}}
Kotzig \cite{Ko} computed the minimum possible number of $C_4$'s in regular tournaments.
We require an analogous lower bound for bipartite Eulerian tournaments. The following lemma provides
this bound.
\begin{lemma}
\label{l21}
Every bipartite Eulerian tournament with vertex class sizes $m$ and $n$
contains at least $m^2n^2/32$ copies of $C_4$.
\end{lemma}
{\bf Proof:}\,
The result is trivial for $m=n=2$. We therefore assume $m+n \geq 6$.
There are four possible ways to orient an undirected $4$-cycle. One of them is $C_4$ and we denote the
other three by $H_1$, $H_2$ and $H_3$ where $H_i$ is the unique orientation having a maximal path
of length $i$. Note that $H_1$ has two sources and two sinks while $H_2$ and $H_3$ each have one
source and one sink. Let $G=(A \cup B,V)$ be a bipartite Eulerian tournament with vertex classes $A$ and $B$
where $|A|=m$ and $|B|=n$.
Let $h_i$ denote the number of distinct copies of $H_i$ in $G$ and let $x$ denote the number
of copies of $C_4$ in $G$. Clearly,
$$
x+h_1+h_2+h_3= {n \choose 2}{m \choose 2}.
$$
Let $t$ denote the total number of sources in all these ${n \choose 2}{m \choose 2}$ graphs.
As each vertex of  $A$ is a source in precisely $(m-1){{n/2} \choose 2}$ such graphs,
and as each vertex of  $B$ is a source in precisely $(n-1){{m/2} \choose 2}$ such graphs,
we have
$$
t=n(n-1){{\frac{m}{2}} \choose 2} + m(m-1){{\frac{n}{2}} \choose 2}.
$$
Every copy of $H_1$ contains two sources, and every copy of $H_2$ and of $H_3$ contains one source.
Thus, $2h_1+h_2+h_3=t$.
It follows that
$$
x-h_1=\frac{mn}{4}(\frac{m}{2}+\frac{n}{2}-1).
$$
For each pair of distinct vertices $u,v$ belonging to the same vertex class, let $k_{u,v}$ denote
the number of common out neighbors. Notice that $k_{u,v}$ is also the number of
common incoming neighbors. Furthermore, if $u,v \in A$ then $n/2-k_{u,v}$ is the number of out neighbors
of $u$ which are incoming neighbors of $v$ and also the number of incoming neighbors
of $u$ which are out neighbors of $v$. The number of $C_4$ copies containing both $u$ and $v$ is
$(n/2-k_{u,v})^2$. The number of $H_1$ copies containing both $u$ and $v$ is $2{{k_{u,v}} \choose 2}$.
An analogous argument holds for $u,v \in B$. As each $H_1$ and each $C_4$ account for two pairs, we have
$$
2x+2h_1=
\left(\sum_{u,v \in A} (\frac{n}{2}-k_{u,v})^2+2{{k_{u,v}} \choose 2}\right) +
\left(\sum_{u,v \in B} (\frac{m}{2}-k_{u,v})^2+2{{k_{u,v}} \choose 2}\right) .
$$
Minimizing the right hand side of the last equality we obtain
$$
2x+2h_1 \ge {m \choose 2}\frac{n^2-2n-1}{8}+{n \choose 2}\frac{m^2-2m-1}{8}.
$$
Therefore,
$$
4x - \frac{mn}{2}(\frac{m}{2}+\frac{n}{2}-1)  \ge
{m \choose 2}\frac{n^2-2n-1}{8}+{n \choose 2}\frac{m^2-2m-1}{8}.
$$
Hence,
$$
x \geq {m \choose 2}\frac{n^2-2n-1}{32}+{n \choose 2}\frac{m^2-2m-1}{32}+
\frac{mn}{8}(\frac{m}{2}+\frac{n}{2}-1)
$$
$$
= \frac{m^2n^2}{32}+\frac{1}{64}(mn^2+nm^2)-\frac{1}{16}mn-\frac{1}{64}(m^2+n^2)+
\frac{1}{64}(m+n) \geq \frac{m^2n^2}{32}.
$$ \endpf

\noindent
We note that Lemma \ref{l21} is asymptotically tight. This can easily be shown by considering random
Eulerian orientations of $K_{m,n}$. We omit the details.

Our next lemma shows that if some arc of a bipartite Eulerian tournament is on too many copies of $C_4$
or on too few copies of $C_4$ then there must be many copies of $C_4$ altogether.
Before we state the lemma we need a few definitions.
Let $G$ be a bipartite Eulerian tournament with vertex class sizes $m$ and $n$.
Let $e=(x,y)$ be an arc of $G$. Let $N^-(x)$ denote the in-neighborhood of $x$
and let $N^+(y)$ denote the out-neighborhood of $y$.
Let $d(e)$ be the number of arcs from $N^+(y)$ to $N^-(x)$.
Notice that $e$ is on precisely $d(e)$ copies of $C_4$ and notice also that there are
$mn/4-d(e)$ arcs from $N^-(x)$ to $N^+(y)$. Let $\alpha(e)=\min\{d(e)/mn~,~1/4-d(e)/mn\}$.
Let $\alpha(G)=\min_{e \in G} \alpha(e)$.

\begin{lemma}
\label{l22}
Let $G$ be a bipartite Eulerian tournament with vertex class sizes $m$ and $n$.
$G$ has at least $m^2n^2(1/16-\alpha(G)+4\alpha(G)^2)$ copies of $C_4$.
\end{lemma}
{\bf Proof:}\,
Let $A$ and $B$ denote the vertex classes of $G$ with $|A|=m$ and $|B|=n$.
Let $e=(x,y)$ be an arc for which $\alpha(e)=\alpha(G)$.
For simplicity let $\alpha=\alpha(e)=\alpha(G)$.
Thus, either $\alpha=d(e)/mn$ or $\alpha=1/4-d(e)/mn$. Without loss of generality
assume that $\alpha=1/4-d(e)/mn$ (the proof of the other case is analogous).
Also without loss of generality assume that $x \in A$ and $y \in B$ (the proof of the other case
is analogous). By our assumptions, $|N^-(x)|=|N^+(x)|=n/2$ and $|N^-(y)|=|N^+(y)|=m/2$.
Also, $e$ appears on precisely $d(e)=mn(1/4-\alpha)$ copies of $C_4$.

Let $T$ be the set of arcs from $N^-(x)$ to $N^+(y)$. Notice that $|T|=\alpha mn$.
Let $S$ be the set of arcs from $N^+(x)$ to $N^-(y)$, let $W$ be the set of arcs from $N^-(y)$ to $N^-(x)$
and let $Z$ be the set of arcs from $N^+(y)$ to $N^+(x)$. Since $G$ is Eulerian we must have
$\alpha mn=|T|=|S|=|W|=|Z|$.

Consider an induced subgraph of $G$ with four vertices, precisely one vertex from
each of $N^+(x), N^-(x), N^+(y), N^-(y)$, and with no arc from $S \cup T \cup W \cup Z$. Such a subgraph
must be a $C_4$. Denote the number of such subgraphs by $x$. Denote the number of $C_4$'s of $G$ having
at least one arc in $S \cup T \cup W \cup Z$ by $y$. Denote the number of $C_4$'s of $G$ with {\em all}
their arcs in $S \cup T \cup W \cup Z$ by $z$. Clearly, $z \leq y$ and there are at least $x+y$
copies of $C_4$ in $G$.

Let $F$ be the set of induced subgraphs of $G$ with four vertices,
precisely one vertex from each of $N^+(x), N^-(x), N^+(y), N^-(y)$. Notice that $|F|=m^2n^2/16$.
Let $F' \subset F$ be the set of elements of $F$ containing at least one arc from
$S \cup T \cup W \cup Z$. Notice that all elements in $F \setminus F'$ are $C_4$'s.
The sum of the arcs of $S \cup T \cup W \cup Z$ over all elements of $F'$ is
precisely $(4\alpha mn)(mn/4)=\alpha m^2n^2$. However, some elements of $F'$ contain more than one
arc from $S \cup T \cup W \cup Z$, and are therefore counted more than once.
In fact, each $e \in T$ and each $e' \in S$ appear
together in precisely one element of $F'$. Similarly, each $e \in W$ and each $e' \in Z$
appear together in precisely one element of $F'$. Thus, there are
$\alpha^2m^2n^2$ elements of $F'$ containing a pair of arcs from $S \cup T$ and there are
$\alpha^2m^2n^2$ elements of $F'$ containing a pair of arcs from $W \cup Z$.
Such elements are counted at least twice. In fact, precisely $z$ of them are counted four times.
It follows that $|F'| \leq \alpha m^2n^2-2\alpha^2m^2n^2-z$. Therefore,
$$
x=|F|-|F'|\geq \frac{1}{16}m^2n^2-\alpha m^2n^2+2\alpha^2m^2n^2+z.
$$

By the definition of $\alpha=\alpha(G)$, each arc of $G$ appears on at least $\alpha mn$ copies of $C_4$.
In particular the sum of the arcs from $S \cup T \cup W \cup Z$ appearing on copies of $C_4$ is
at least $(4\alpha mn)(\alpha mn)=4\alpha^2 m^2n^2$. Precisely $z$ of the $C_4$'s containing arcs
from $S \cup T \cup W \cup Z$ are counted $4$ times in this way. Notice that if some $C_4$ has
three arcs from $S \cup T \cup W \cup Z$ then it must have four arcs in $S \cup T \cup W \cup Z$.
Thus, $y \geq z+(4\alpha^2 m^2n^2-4z)/2=2\alpha^2 m^2n^2-z$. Consequently,
$$
x+y \geq m^2n^2(\frac{1}{16}-\alpha+4\alpha^2).
$$
\endpf

Recall that a $k$-uniform hypergraph is a family of $k$-subsets of some set $V$ of vertices.
The $k$-subsets are the {\em edges} of the hypergraph.
If $x,y$ are two vertices of a hypergraph then let $d(x)$ denote the number of edges containing $x$
and let $d(x,y)$ denote the number of edges containing both $x$ and $y$ (their {\em co-degree}).
A proper $r$-coloring of the edges of a hypergraph is a partition of the set of edges to $r$ subsets
such that any two edges in the same subset are disjoint.
Kahn \cite{Ka} has proved a very powerful result giving an upper bound for the minimum number of
colors in a proper edge-coloring of a uniform hypergraph (his result is, in fact, more general).
\begin{lemma}[Kahn \cite{Ka}]
\label{l23}
For every $k \geq 2$ and every $\epsilon > 0$ there exists a $\delta > 0$ such that the following
statement is true:\\
If $H$ is a $k$-uniform hypergraph on a set $V$ of vertices, and if
$d(v) \leq D$ for all $v \in V$ and $d(u,v) < \delta D$ for all distinct $u,v \in V$,
then there is a proper coloring of the edges of $H$ with at most $(1+\epsilon)D$ colors. \endpf
\end{lemma}

\noindent
{\bf Completing the proof of Theorem \ref{t1}.}
Let $\epsilon > 0$. We need to prove that if $m$ is sufficiently large as a function of $\epsilon$ and
if $n \geq m$ and $G$ is a bipartite Eulerian tournament with vertex class sizes $m$ and $n$,
then $G$ has at least $\frac{1}{4+\sqrt 8}mn(1-\epsilon)$ arc-disjoint copies of $C_4$.

We construct a $4$-uniform hypergraph $H$ as follows. The vertices of $H$ are the $mn$ arcs
of $G$ and the edges of $H$ are the copies of $C_4$ in $G$. Let $e=(x,y)$ be an arc of $G$
such that $d(e)$ is maximal. Clearly, $d(e) \leq mn/4$. By Lemma \ref{l21} we also have
$d(e) \geq 4(m^2n^2/32)/(mn) \geq mn/8$.
In particular $\alpha(e)=1/4-d(e)/mn$ and $0 \leq \alpha(G) \leq \alpha(e) \leq 1/8$.
Notice that any two arcs of $G$ appear together on at most $\max \{m/2, n/2\}=n/2$ copies of $C_4$.
Thus, with $D=d(e)$ and with $\delta=\delta(4,\epsilon)$ chosen as in Lemma \ref{l23} we
have that for $m$ sufficiently large, the hypergraph $H$ satisfies the conditions of Lemma \ref{l23}.
In particular, there is a coloring of the $C_4$ copies of $G$ with at most $(1+\epsilon)d(e)$
colors such that each color class corresponds to a set of arc-disjoint copies of $C_4$.
By Lemma \ref{l21} and Lemma \ref{l22}, the total number of edges of $H$ is at least
$$
\max\{ \frac{1}{32}m^2n^2 ~,~ m^2n^2(\frac{1}{16}-\alpha(G)+4\alpha(G)^2) \}.
$$
Thus, there must be a set of arc-disjoint copies of $C_4$ with at least
$$
\max\{ \frac{m^2n^2/32}{(1+\epsilon)d(e)} ~,~ 
\frac{m^2n^2(1/16-\alpha(G)+4\alpha(G)^2)}{(1+\epsilon)d(e)} \}
$$
elements. Since $\alpha(G) \leq \alpha \leq 1/8$ it suffices to prove that
$$
\frac{mn}{1+\epsilon}\max\{ \frac{1/32}{1/4-\alpha(e)} ~,~ 
\frac{1/16-\alpha(e)+4\alpha(e)^2}{1/4-\alpha(e)} \} \geq \frac{1}{4+\sqrt 8}mn(1-\epsilon).
$$
Indeed, this follows from the fact that for $0 \leq z \leq 1/8$
the value of $\max\{(1/32)/(1/4-z)~,~(1/16-z+4z^2)/(1/4-z)\}$
is minimized when $z=1/8-\sqrt{1/128}$.
\endpf

\section{Proof of Theorem \ref{t2}}
The proofs of Lemma \ref{l21} and Lemma \ref{l22}, and hence of Theorem \ref{t1}, can easily
be modified to hold for {\em almost} Eulerian bipartite tournaments. More precisely,
we say that a bipartite orientation of $K_{m,n}$ is $\delta$-Eulerian if the outdegree and indegree of
each vertex from the $m$-vertex class is at least $(1-\delta)n/2$ and if the outdegree and indegree of
each vertex from the $n$-vertex class is at least $(1-\delta)m/2$.
Given $\epsilon > 0$, it is straightforward to check that Lemma \ref{l21} holds for $\delta$-Eulerian
bipartite tournaments if one replaces the $m^2n^2/32$ bound with $m^2n^2(1-\epsilon)$ provided that
$\delta=\delta(\epsilon)$ is sufficiently small and $\min \{m,n\}$ is sufficiently large as a function
of $\epsilon$. It is straightforward to check that a similar relaxation also holds for Lemma \ref{l22}.
Thus, we obtain the following relaxation for Theorem \ref{t1}.
\begin{lemma}
\label{l31}
Let $\epsilon > 0$. There exists $\delta=\delta(\epsilon) > 0$ and $M = M(\epsilon)$ such that
if $M \leq m \leq n$, any $\delta$-Eulerian orientation of $K_{m,n}$ contains at least
$\frac{1}{4+\sqrt 8}mn(1-\epsilon)$ arc-disjoint copies of $C_4$.
\end{lemma}

\noindent
{\bf Completing the proof of Theorem \ref{t1}.}
Let $\epsilon > 0$. We need to show that there exists $N=N(\epsilon)$ such that for all $n > N$,
any regular tournament with $n$ vertices has at least $\frac{1}{8+\sqrt{32}}n^2(1-\epsilon)$
arc-disjoint directed $4$-cycles.

We choose $\delta=\delta(\epsilon/2)$ to be as in Lemma \ref{l31}. Choose
$N \geq \max \{2M^2~,~64/\epsilon^2\}$, where
$M=M(\epsilon/2)$ is the constant from Lemma \ref{l31}, to be the smallest integer satisfying
$4N^2\exp(-\delta^2 N/128) < 1$.
Let $n > N$. The asymptotic nature of our proof allows us to assume, for simplicity, that $n$
is a square of an integer. Therefore, let $m=\sqrt{n}$. Let $T$ be a regular tournament with $n$ vertices.
Each $v \in T$ chooses uniformly and independently an integer from $\{1,\ldots,m\}$.
Let $V_i$ denote the set of vertices that chose $i$, and for $1 \leq i < j \leq m$ let $G_{i,j}$ be the
complete bipartite tournament with vertex classes $V_i$ and $V_j$.

We will prove that with positive probability, {\em all} the ${m \choose 2}$ bipartite tournaments
$G_{i,j}$ are $\delta$-Eulerian, and that {\em all} $V_i$ have size at least $M$ and at most $2m$.
This suffices to complete the theorem since by Lemma \ref{l31} we have that the total number
of arc-disjoint $C_4$ in $T$ is at least
$$
\sum_{1 \leq i < j \leq m} \frac{1}{4+\sqrt 8}|V_i||V_j|(1-\frac{\epsilon}{2})
$$
$$
=  \frac{1}{8+\sqrt {32}}(1-\frac{\epsilon}{2})(n^2-\sum_{i=1}^m |V_i|^2)
\geq \frac{1}{8+\sqrt {32}}(1-\frac{\epsilon}{2})(n^2-4m^3) \geq
$$
$$
\geq \frac{1}{8+\sqrt {32}}(1-\frac{\epsilon}{2})n^2(1-\frac{\epsilon}{2})
\geq \frac{1}{8+\sqrt {32}}n^2(1-\epsilon).
$$

Fix a vertex $v$. Let $d^+_i(v)$ (resp. $d^-_i(v)$) denote the number of out (resp. incoming) neighbors
of $v$ in $V_i$.
The expectation of $d^+_i(v)$ (resp. $d^-_i(v)$) is precisely $m^{-1}(n-1)/2=(n-1)/(2m)$.
A large deviation inequality of Chernoff (cf. \cite{AlSp} Appendix A) bounds the probability that a random
variable which is the sum of independent indicator random variables deviates from its expectation by more
than a specified amount. More precisely, for the variable $d^+_i(v)$ and the amount
$\frac{\delta}{4}\frac{n-1}{2m}$ we get
$$
\Pr\left[\left| d^+_i(v) - \frac{n-1}{2m} \right| > \frac{\delta}{4}\frac{n-1}{2m}\right] <
2 \exp \left(
-\frac{\delta^2(n-1)^2/64m^2}{2(n-1)/2m}+
\frac{\delta^3(n-1)^3/512m^3}{2(n-1)^2 / 4m^2}\right) =
$$
$$
2 \exp \left(-\frac{\delta^2(n-1)}{64m}+\frac{\delta^3(n-1)}{128m} \right) < \frac{1}{2mn}.
$$
Precisely the same bound holds for $d^-_i(v)$. Since there are $n$ choices for vertices and $m$ choices
for $i$, we get that with probability {\em greater than} $1-2 \cdot n \cdot m \cdot \frac{1}{2mn} = 0$,
for {\em all} vertices $v$ and all indices $i$,
$$
\left| d^+_i(v) - \frac{n-1}{2m} \right| \leq \frac{\delta}{4}\frac{n-1}{2m} \qquad
\left| d^-_i(v) - \frac{n-1}{2m} \right| \leq \frac{\delta}{4}\frac{n-1}{2m}.
$$
We therefore {\em fix} a partition $V_1,\ldots,V_m$ for which all of these conditions are met.
In particular, this means that for all $i=1,\ldots,m$,
$$
2m > (1+\frac{\delta}{4})\frac{n-1}{m} \geq |V_i| \geq (1-\frac{\delta}{4})\frac{n-1}{m} > M,
$$
and that for all vertices $v$ and all $i=1,\ldots,m$
$$
\min \{d^-_i(v)~,~d^+_i(v)\} \geq (1-\frac{\delta}{4})\frac{n-1}{2m} \geq (1-\delta)(1+\frac{\delta}{4})\frac{n-1}{2m} \geq
\frac{|V_i|}{2}(1-\delta).
$$
Thus, all the ${m \choose 2}$ bipartite tournaments
$G_{i,j}$ are $\delta$-Eulerian, and all $V_i$ have size at least $M$ and at most $2m$.
\endpf

\end{document}